\nonstopmode
\documentclass[a4paper, 10pt, reqno]{amsart}
\usepackage{latexsym}
\usepackage{fancyhdr}
\usepackage{amsmath, amssymb}
\usepackage[ansinew]{inputenc}
\usepackage[all]{xy}
\usepackage{verbatim}
\usepackage{hyperref}
\swapnumbers
\theoremstyle{plain}
\newtheorem*{lemma*}{Lemma}

\newtheorem*{theorem*}{Theorem}

\newtheorem*{proposition*}{Proposition}

\newtheorem*{corollary*}{Corollary}

\newtheorem*{claim*}{Claim}

\theoremstyle{definition}
\newtheorem*{definition*}{Definition}

\newtheorem*{example*}{Example}

\newtheorem*{algorithm*}{Algorithm}
\newtheorem*{remark*}{Remark}
\newtheorem*{remarks*}{Remarks}

\numberwithin{equation}{subsection}
\newenvironment{demo}[1]{\par\smallskip\noindent{\bf #1.}}{\par\smallskip}

\newenvironment{proclaim}[1]{\par\medskip\noindent{\bf #1.}\it}{\par\smallskip}

\numberwithin{equation}{section}

\def\al{\alpha}

\def\ga{\gamma}

\def\la{\lambda}

\def\C{\mathbb{C}}

\def\R{\mathbb{R}}
\def\RR{\mathbb{R}}

\def\x{\times}

\def\oo{\infty}

\def\<{\langle}
\def\>{\rangle}

\def\sr#1%
{\ifmmode{}^\dagger\else${}^\dagger$\fi\ifvmode
\vbox to 0pt{\vss
 \hbox to 0pt{\hskip\hsize\hskip1em
 \vbox{\hsize3cm\eightpoint\raggedright\pretolerance10000
 \noindent #1\hfill}\hss}\vss}\else
 \vadjust{\vbox to0pt{\vss%
 \hbox to 0pt{\hskip\hsize\hskip1em%
 \vbox{\hsize3cm\eightpoint\raggedright\pretolerance10000%
 \noindent #1\hfill}\hss}\vss}}\fi%
}

\begin{document}

\title[Many parameter H\"older perturbation of unbounded operators]
{Many parameter H\"older perturbation of unbounded operators
}
\author{Andreas Kriegl, Peter W. Michor, Armin Rainer}
\address{
Andreas Kriegl: 
Fakult\"at f\"ur Mathematik, Universit\"at Wien,
Nordbergstrasse 15, A-1090 Wien, Austria
}
\email{Andreas.Kriegl@univie.ac.at}
\address{
Peter W. Michor:
Fakult\"at f\"ur Mathematik, Universit\"at Wien,
Nordbergstrasse 15, A-1090 Wien, Austria
}
\email{Peter.Michor@univie.ac.at}
\address{
Armin Rainer: Department of Mathematics, University of Toronto, 
40 St.\ George Street, Toronto, Ontario, Canada M5S 2E4
}
\email{armin.rainer@univie.ac.at}
\date{\today}
\thanks{P.W.M. was supported by FWF grant P~21030-N13; A.R. by  FWF grant J~2771}

\begin{abstract}
If $u\mapsto A(u)$ is a $C^{0,\al}$-mapping, for $0< \al \le 1$, having as values unbounded
self-adjoint
operators with compact resolvents and common domain
of definition, parametrized by $u$ in an (even infinite dimensional)
space, then any continuous (in $u$) arrangement of the eigenvalues of 
$A(u)$ is indeed $C^{0,\al}$ in $u$.
\end{abstract}
\maketitle

\begin{theorem*}
Let $U \subseteq E$ be a $c^\infty$-open subset in a convenient vector
space $E$, and $0 < \al \le 1$.
Let $u \mapsto A(u)$, for $u \in U$, be a $C^{0,\al}$-mapping with values unbounded
self-adjoint operators
in a Hilbert space $H$ with common domain of definition and with compact
resolvent.
Then any (in $u$) continuous eigenvalue $\la(u)$ of $A(u)$ 
is $C^{0,\al}$ in $u$.
\end{theorem*}

\subsection*{Remarks and definitions}
This paper is a complement to \cite{KM03} and builds upon it. 
A function $f:\mathbb R\to\mathbb R$ is called $C^{0,\al}$ if 
$\frac{f(t)-f(s)}{|t-s|^\al}$ is
locally bounded in $t\ne s$. 
For $\al=1$ this is Lipschitz.

Due to \cite{Boman67} a mapping $f : \R^n \to \R$ is $C^{0,\al}$ if and
only if $f\circ c$ is $C^{0,\al}$ for each
smooth (i.e.\ $C^\oo$) curve $c$.
\cite{Faure89} has shown that this holds for even more general concepts of
H\"older differentiable
maps.

%Due to \cite{FK88}, see also \cite[12.7]{KM97} 
%a mapping between Banach spaces is locally Lipschitz if and
%only if $f\circ c$ is  Lipschitz for each
%smooth (i.e.\ $C^\oo$) curve $c$.

A convenient vector space (see \cite{KM97}) is a locally convex vector space
$E$ satisfying the following equivalent conditions: Mackey Cauchy sequences
converge; $C^\infty$-curves in $E$ are locally integrable in $E$; 
a curve $c:\mathbb R\to E$ is $C^\infty$ (Lipschitz) if and only if 
$\ell\circ c$ is $C^\infty$ (Lipschitz)
for all continuous linear functionals $\ell$.
The $c^\infty$-topology on $E$
is the final topology with respect to all smooth curves (Lipschitz curves).
Mappings $f$ defined on open (or even $c^\infty$-open) 
subsets of convenient vector spaces $E$ 
are called $C^{0,\al}$ (Lipschitz) if $f\circ c$ is $C^{0,\al}$ (Lipschitz) for every smooth curve
$c$. If $E$ is a Banach space then a $C^{0,\al}$-mapping is locally H\"older-continuous
of order $\al$ in the usual sense. This has been proved in \cite{Faure91},
which is not easily accessible, thus we include a proof in the lemma below. 
For the Lipschitz case see \cite{FK88} and \cite[12.7]{KM97}.

That a mapping $t\mapsto A(t)$ defined on a $c^\oo$-open subset $U$ of a
convenient vector space $E$
is  %real analytic, $C^\infty$, or 
$C^{0,\al}$ with values in
unbounded operators means the following:
There is a dense subspace $V$ of the Hilbert space $H$
such that $V$ is the domain of definition of each $A(t)$, and such
that $A(t)^*=A(t)$. And furthermore,
$t\mapsto \langle A(t)u,v\rangle$
is %real analytic, $C^\infty$, or 
$C^{0,\al}$ for each $u\in V$ and $v\in H$
in the sense of the definition given above.

This implies that $t\mapsto A(t)u$ is of the same class $U\to H$
for each $u\in V$ by \cite[2.3]{KM97}, \cite[2.6.2]{FK88}, or\cite[4.1.14]{Faure91}.
This is true because $C^{0,\al}$ can be described by
boundedness conditions only; and for these the uniform boundedness
principle is valid.

\begin{lemma*} {\rm \cite{Faure91}}
Let $E$ and $F$ be Banach spaces, $U$ 
open in $E$. Then, a mapping $f:U\to F$ is $C^{0,\al}$ if and only if $f$ 
is locally H\"older of order $\al$, i.e., 
$\frac{\|f(x)-f(y)\|}{\|x-y\|^\al}$ is locally 
bounded.
\end{lemma*}

\begin{demo}{Proof}
If $f$ is $C^{0,\al}$ but not locally H\"older near 
$z\in U$, there are $x_n \ne y_n$ in $U$ with 
$\|x_n-z\|\le 1/{4^n}$ and $\|y_n-z\|\le 1/{4^n}$, such that 
$\|f(y_n)-f(x_n)\|\ge n.2^n.\|y_n- x_n\|^\al$. Now we apply the general 
curve lemma \cite[12.2]{KM97} with $s_n:=2^n.\|y_n- x_n\|$ and 
$c_n(t):= x_n-z+t\frac{y_n- x_n}{2^n\|y_n-x_n\|}$ to get a smooth curve $c$ 
with $c(t+t_n)-z=c_n(t)$ for $0\le t \le s_n$. Then
$\frac1{s_n^\al}\|(f\circ c)(t_n+s_n)-(f\circ c)(t_n)\|
     =\frac1{2^{n\al}.\|y_n-x_n\|^\al} \|f(y_n)-f(x_n)\| \ge n$.
The converse is obvious.
\qed
\end{demo}

\subsection*{The theorem holds for $E=\RR$}
Let $t \mapsto A(t)$ be a $C^{0,\al}$-curve.
Going through the proof of the resolvent lemma in \cite{KM03} carefully,
we find that $t \mapsto A(t)$ is a
$C^{0,\al}$-mapping $U \to L(V,H)$, and thus the resolvent $(A(t)-z)^{-1}$ is
$C^{0,\al}$ into $L(H,H)$ in
$t$ and $z$ jointly.

For a continuous eigenvalue $t\mapsto \la(t)$ as in the theorem,
let the eigenvalue $\la(s)$ of $A(s)$ have multiplicity $N$ for $s$ fixed.
Choose a simple closed curve $\ga$ in the resolvent set of $A(s)$ enclosing
only $\la(s)$ among all eigenvalues of $A(s)$.
Since the global resolvent set $\{(t,z) \in \R \times \C : (A(t)-z) : V \to
H ~\text{is invertible}\}$ is open,
no eigenvalue of $A(t)$ lies on $\ga$, for $t$ near $s$. Consider
\[
t \mapsto -\frac{1}{2 \pi i} \int_\ga (A(t)-z)^{-1} dz =:P(t),
\]
a $C^{0,\al}$-curve of projections (on the direct sum of all eigenspaces
corresponding to eigenvalues in the interior
of $\ga$) with finite dimensional ranges and constant ranks.
So for $t$ near $s$, there are equally many eigenvalues 
(repeated with multiplicity) in the interior of
$\ga$. 
Let us order them by size, $\mu_1(t) \le \mu_2(t) \le \cdots \le \mu_N(t)$,
for all $t$. 
The image of $t \mapsto P(t)$, for $t$ near $s$ describes a finite
dimensional $C^{0,\al}$ vector subbundle of
$\R \times H \to \R$, since its rank is constant.
The set $\{\mu_i(t) : 1 \le i \le N\}$ 
represents the eigenvalues of $P(t)A(t)|_{P(t)(H)}$.
By the following result, it forms a $C^{0,\al}$-parametrization 
of the eigenvalues of $A(t)$ inside $\ga$, for $t$ near $s$.

The eigenvalue $\la(t)$ is a continuous (in $t$) choice among the $\mu_i(t)$, 
and it is $C^{0,\al}$ in $t$ by the proposition below.

\begin{proclaim}{Result}
\emph{(\cite{Weyl12}, see also \cite[III.2.6]{Bhatia97})}
Let $A,B$ be $N \times N$ Hermitian matrices. 
Let  $\mu_1(A)  \le  \mu_2(A)  \le  \cdots  \le  \mu_N(A)$  and
$\mu_1(B)  \le  \mu_2(B)  \le  \cdots  \le \mu_N(B)$ denote the
eigenvalues of $A$ and $B$, respectively.
Then 
\[
\max_j |\mu_j(A)-\mu_j(B)| \le \|A-B\|.
\]
Here $\|.\|$ is the operator norm.
\end{proclaim}

\begin{proposition*}
Let $0 < \al \le 1$.
Let $U\ni u\mapsto A(u)$ be a $C^{0,\al}$-mapping of Hermitian $N\x N$ matrices.
Let $u\mapsto\la_i(u)$, $i=1,\dots,N$ be continuous mappings 
which together parametrize the eigenvalues of $A(u)$. 
Then each $\la_i$ is $C^{0,\al}$. 
\end{proposition*}

\begin{demo}{Proof}
It suffices to check that $\la_i$ is $C^{0,\al}$ along each smooth curve in $U$,
so we may assume without loss that $U=\mathbb R$.
We have to show that each continuous eigenvalue $t\mapsto \la(t)$ 
is a $C^{0,\al}$-function on each compact interval $I$ in $U$. 
Let $\mu_1(u)\le \dots\le \mu_N(u)$ be the increasingly ordered arrangement of 
eigenvalues. 
Then each $\mu_i$ is a $C^{0,\al}$-function on $I$ with a common H\"older  
constant $C$ by the result above.
Let $t<s$ be in $I$.
Then there is an $i_0$ such that $\la(t)=\mu_{i_0}(t)$. Now let $t_1$ be the 
maximum of all $r\in [t,s]$ such that $\la(r)=\mu_{i_0}(r)$. If $t_1<s$ then 
$\mu_{i_0}(t_1)=\mu_{i_1}(t_1)$ for some $i_1\ne i_0$. 
Let $t_2$ be the 
maximum of all $r\in [t_1,s]$ such that $\la(r)=\mu_{i_1}(r)$. If $t_2<s$ then 
$\mu_{i_1}(t_2)=\mu_{i_2}(t_2)$ for some $i_2\notin \{i_0,i_1\}$. 
And so on until $s=t_k$ for some $k\le N$. Then we have (where $t_0=t$)
\begin{align*}
\frac{|\la(s)-\la(t)|}{(s-t)^\al} 
&\le \sum_{j=0}^{k-1} \frac{|\mu_{i_j}(t_{j+1})-\mu_{i_j}(t_j)|}{(t_{j+1}-t_j)^\al}
\cdot \left(\frac{t_{j+1}-t_j}{s-t}\right)^\al
\le C k
\le C N.
\qed
\end{align*}
\end{demo}

\subsection*{Proof of the theorem}
For each smooth curve $c:\mathbb R\to U$ the curve
$\mathbb R\ni t\mapsto A(c(t))$ is $C^{0,\al}$,
and by the
1-parameter case the eigenvalue $\la(c(t))$ is $C^{0,\al}$.
But then $u\mapsto \la(u)$ is $C^{0,\al}$.
\qed

\subsection*{Remark} 
Let $u \mapsto A(u)$ be Lipschitz. 
Choose a fixed continuous ordering of the roots, e.g., by size.
We claim that along a smooth or Lipschitz curve $c(t)$ in $U$, none of these can accelerate 
to $\infty$ of $-\infty$ in finite time. 
Thus we may denote them as
$\dots\la_i(u)\le \la_{i+1}(u)\le\dots$, for all $u\in U$.
Then each $\la_i$ is Lipschitz.

The claim can be proved as follows:
Let $t \mapsto A(t)$ be a Lipschitz curve.
By reducing to the projection $P(t)A(t)|_{P(t)(H)}$, we may assume that $t \mapsto A(t)$ is  a Lipschitz curve of $N \times N$ Hermitian matrices. 
So $A'(t)$ exists a.e.\ and is locally bounded.
Let $t \mapsto \la(t)$ be a continuous eigenvalue.
It follows that $\la$ satisfies \cite[(6)]{KM03} a.e.\ and, as in the proof of \cite[(7)]{KM03}, one shows that for each compact interval $I$ there is a 
constant $C$ such that $|\la'(t)| \le C+C |\la(t)|$ a.e.\ in $I$.
Since $t \mapsto \la(t)$ is Lipschitz, in particular, absolutely continuous, Gronwall's lemma (e.g.\ \cite[(10.5.1.3)]{Dieudonne60})  implies that 
$|\la(s) - \la(t)| \le (1+|\la(t)|) (e^{a |s-t|}-1)$ for a constant $a$ depending only on $I$.

%\nocite{FF89}
%\bibliography{biblio}
%\bibliographystyle{amsplain}

\def\cprime{$'$}
\providecommand{\bysame}{\leavevmode\hbox to3em{\hrulefill}\thinspace}
\providecommand{\MR}{\relax\ifhmode\unskip\space\fi MR }
\providecommand{\MRhref}[2]{%
  \href{http://www.ams.org/mathscinet-getitem?mr=#1}{#2}
}
\providecommand{\href}[2]{#2}

\end{document}